\documentclass[11pt,a4paper]{article} %replace-string
\usepackage{amsmath}
\usepackage{amsfonts}
\usepackage{amssymb}
\usepackage{amsthm}
\usepackage{amscd}
\usepackage{graphicx}
\theoremstyle{plain}
\newtheorem{lemma}{Lemma}[section]

\newtheorem{theorem}{Theorem}[section]
\theoremstyle{definition}

\newtheorem{step}{STEP}

\title{Global properties of indefinite metrics with parallel Weyl tensor}

\author{
	\emph{Andrzej Derdzinski}
	Ohio State University \\
	Department of Mathematics \\ 
	Columbus, OH 43210, USA \\ 
	e-mail: \emph{andrzej@math.ohio-state.edu} \\
	url: \emph{
             http:/\hskip-1.5pt/www.math.ohio-state.edu/\~{}\hskip-1ptandrzej}
\and
	\emph{Witold Roter}
	Wroc\l aw University of Technology \\
	Institute of Mathematics and Computer Science  \\ 
	Wy\-brze\-\.ze Wys\-pia\'n\-skiego 27, 50-370 Wroc\l aw, Poland \\ 
	e-mail: \emph{Witold.Roter@pwr.wroc.pl} \\
	url: \emph{http:/\hskip-1.5pt/www.im.pwr.wroc.pl/ang/visit/wroter.html}
       }

\date{}

\begin{document}

\maketitle

\begin{abstract}
This is an exposition of some recent results on ECS manifolds, by which 
we mean pseudo-Riemannian manifolds of dimensions greater than 3 that are 
neither conformally flat nor locally symmetric, and have parallel Weyl 
tensor. All ECS metrics are indefinite. We state two classification theorems, 
describing the local structure of ECS manifolds, and outline an argument 
showing that compact ECS manifolds exist in infinitely many dimensions greater 
than  4. We also discuss some properties of compact manifolds that admit ECS 
metrics, and provide a list of open questions about compact ECS manifolds.
\end{abstract}

\section*{Introduction}%\label{derdzinski-roter:intr}
An {\it ECS manifold\/} is any pseu\-\hbox{do\hskip.7pt-}Riem\-ann\-i\-an 
manifold of dimension $\,n\ge4\,$ that has parallel Weyl tensor 
($\nabla W\hskip-1.7pt=0$) without being con\-for\-mal\-ly flat 
($W\hskip-1.7pt=0$) or locally symmetric ($\nabla\hskip-1.5ptR=0$). ECS 
manifolds exist in every dimension $\,n\ge4$, and their metrics are all 
indefinite \cite[Corol\-lary~3]{derdzinski-roter:r}, 
\cite[Theorem~2]{derdzinski-roter:dr77}, cf.\ 
\cite[Re\-mark~16.75(iii)]{derdzinski-roter:b}. This paper presents selected 
results on ECS manifolds.

A local classification of ECS metrics, described in 
\S\S\ref{derdzinski-roter:lsdo}\hskip2pt--\hskip2pt\ref{derdzinski-roter:lsdt}, 
is naturally divided into two cases, depending on the dimension of a null 
parallel distribution first introduced by Ol\-szak \cite{derdzinski-roter:o}. 
In \S\ref{derdzinski-roter:cxis} we state a theorem about the existence of 
compact ECS manifolds in infinitely many dimensions $\,n\ge5$. An outline of 
its proof is given in \S\ref{derdzinski-roter:pfth}. 
Sections~\ref{derdzinski-roter:cprp}\hskip2pt--\hskip2pt\ref{derdzinski-roter:comm} 
contain results on general properties of compact manifolds admitting ECS 
metrics, and a list of open questions concerning compact ECS manifolds.

Metrics with parallel Weyl tensor were first studied by Chaki and Gupta 
\cite{derdzinski-roter:cg}, who referred to them as {\it con\-for\-mal\-ly 
symmetric}. Although the latter term was adopted by many authors 
(\cite{derdzinski-roter:dh,derdzinski-roter:h,derdzinski-roter:s,derdzinski-roter:sky}, to name just a few), 
it may be misleading, since it appears in the literature with at least two 
other meanings \cite{derdzinski-roter:hm,derdzinski-roter:n}. This is why, 
rather than speaking of con\-for\-mal symmetry, we use the phrases `parallel 
Weyl tensor' and `ECS manifold' (or, `metric'), `ECS' being abbreviated from 
{\it essentially con\-for\-mal\-ly symmetric}.

We wish to express our gratitude to David Blair and Ryszard Deszcz for 
comments concerning terminology and pseu\-\hbox{do\hskip.7pt-}sym\-me\-try of 
ECS metrics. We also thank Maciej Dunajski for bringing to our attention the 
results of \cite{derdzinski-roter:dw}, and Tadeusz Januszkiewicz for pointing 
out the fact described in the last paragraph of \S\ref{derdzinski-roter:open}.

\section{Parallel Weyl tensor}\label{derdzinski-roter:pwte}
The symbols $\,R\,$ and $\,W$ always denote the curvature tensor and Weyl 
con\-for\-mal tensor of the pseu\-\hbox{do\hskip.7pt-}Riem\-ann\-i\-an metric 
in question, while $\,\rho,\hskip.7pt\text{\rm s}\hskip.7pt\,$ and 
$\,\hskip.4pt\nabla\,$ stand for its Ric\-ci tensor, scalar curvature and 
Le\-vi-Ci\-vi\-ta connection. In dimensions $\,n\ge4$, one has the 
decomposition
\begin{equation}\label{derdzinski-roter:dec}
\begin{array}{ccccccc}
R&=&S&+&E&+&W\hskip50pt\\
\mathrm{curvature}&=&\mathrm{scalar}&+&\mathrm{Ein\-stein}&+&\mathrm{Weyl}
\hskip50pt
\end{array}
\end{equation}
of $\,R\,$ into its irreducible components under the action of the 
pseu\-\hbox{do\hskip.7pt-}or\-thog\-o\-nal group; if $\,n=4$, the underlying 
manifold is oriented, and the metric signature is Riemannian or neutral, $\,W$ 
can be further decomposed into its self-du\-al and anti-self-du\-al parts. The 
simplest linear conditions imposed on $\,R\,$ are $\,S=0\,$ 
(sca\-lar-flat\-ness), $\,W\hskip-1.7pt=0\,$ (con\-for\-mal flatness) and 
$\,E=0\,$ (which defines Ein\-stein metrics). The decomposition of 
$\,\nabla\hskip-1.5ptR\,$ resulting from (\ref{derdzinski-roter:dec}) leads to 
the analogous conditions $\,\nabla\hskip-.8ptS=0$, 
$\hskip.7pt\nabla\hskip-1.3ptE=0\,$ and $\,\nabla W\hskip-1.7pt=0$. Questions 
about metrics with parallel Weyl tensor may, consequently, be considered as 
natural as those concerning constant scalar curvature, or parallel Ric\-ci 
tensor (including the case of Ein\-stein metrics and their products).

The requirement that the Weyl tensor be parallel is also related to some other 
conditions that are of independent interest. For instance, all ECS metrics 
are sca\-lar-flat \cite[Theorem~7]{derdzinski-roter:dr78}, which, combined 
with the second Bianchi identity and the relation $\,\nabla W\hskip-1.7pt=0$, 
implies in turn that they have harmonic curvature \cite{derdzinski-roter:b} 
(in other words, the Ric\-ci tensor satisfies the Codazzi equation).

Furthermore, all ECS manifolds are sem\-i\-sym\-met\-ric 
\cite[Theorem~9]{derdzinski-roter:dr78}, and hence 
pseu\-do\-sym\-met\-ric \cite{derdzinski-roter:dh}. Finally, in terms of the 
invariant $\,d\in\{1,2\}\,$ appearing in relation (\ref{derdzinski-roter:dot}) 
below, every ECS metric has low co\-ho\-mo\-ge\-ne\-i\-ty (at most 
$\,d\hskip.7pt$), an ECS manifold with $\,d=1\,$ is necessarily 
Ric\-ci-re\-cur\-rent, and, in dimension four, ECS metrics with $\,d=2\,$ are 
all self-du\-al. See 
\S\S\ref{derdzinski-roter:odis}\hskip2pt--\hskip2pt\ref{derdzinski-roter:lsdo}.

\section{The Ol\-szak distribution}\label{derdzinski-roter:odis}
Let a pseu\-\hbox{do\hskip.7pt-}Riem\-ann\-i\-an manifold $\,(M,g)\,$ of 
dimension $\,n\ge4\,$ have parallel Weyl tensor $\,W\hskip-1.4pt$. The {\it 
Ol\-szak distribution\/} of $\,(M,g)$, introduced by Ol\-szak 
\cite{derdzinski-roter:o}, is the sub\-bun\-dle $\,\mathcal{D}\,$ of 
$\,T\hskip-.3ptM\,$ such that the sections of $\,\mathcal{D}\,$ are precisely 
those vector fields $\,u\,$ which satisfy the condition
\begin{equation}\label{derdzinski-roter:gww}
g(u,\,\cdot\,)\,\wedge\,W(v,v\hskip.4pt'\hskip-1pt,\,\cdot\,,\,\cdot\,)\,
\hskip.7pt=\hskip.7pt\,0
\end{equation}
for all vector fields $\,v,v\hskip.4pt'$ (exterior multiplication of a 
$\,1$-form by a $\,2$-form). At every point $\,x\in M$, the space 
$\,\mathcal{D}\hskip-.3pt_x$ thus consists, in addition to the zero vectors, 
of all vectors that $\,\wedge$-divide each $\,2$-form in the image of 
$\,W_{\hskip-.7ptx}$ acting on $\,2$-forms. In other words, 
$\,\mathcal{D}\hskip-.3pt_x$ is the subspace of $\,T_xM$ formed by the 
vectors that lie in the image of every nonzero $\,2$-form in the image of 
$\,W_{\hskip-.7ptx}$. (We use $\,g_x$ to identify $\,2$-forms at $\,x\,$ 
with skew-ad\-joint endomorphisms of $\,T_xM$.)

The distribution $\,\mathcal{D}\,$ is obviously parallel. We denote its 
dimension by $\,d\hskip.4pt$. For any $\,n$-di\-men\-sion\-al 
pseu\-\hbox{do\hskip.7pt-}Riem\-ann\-i\-an manifold $\,(M,g)\,$ with parallel 
Weyl tensor, $\,d\in\{0,1,2,n\}$, while $\,d=n\,$ if and only if $\,g\,$ is 
con\-for\-mal\-ly flat. Furthermore,
\begin{equation}\label{derdzinski-roter:dot}
\text{\rm in\ every\ ECS\ manifold,}\hskip7pt\mathcal{D}\hskip5pt\text{\rm is\ 
a\ null\ parallel\ distribution\ of\ dimension}\hskip6pt
d\,\in\,\{1,2\}\hskip.7pt.
\end{equation}
The facts just stated are due to Ol\-szak \cite{derdzinski-roter:o}. Their 
proofs an also be found in \cite[Lemma~2.1]{derdzinski-roter:dr07c}.
\begin{lemma}\label{derdzinski-roter:detwo}For any 
pseu\-\hbox{do\hskip.7pt-}Riem\-ann\-i\-an manifold\/ $\,(M,g)\,$ with 
parallel Weyl tensor, one has\/ $\,d=2\hskip2pt\,$ if and only if\/ 
$\,\hskip1ptW\hskip-1.7pt=\pm\hskip2.5pt\omega\otimes\omega\hskip1pt\,$ for 
some sign\/ $\,\pm\,$ and some parallel differential\/ $\,2$-form\/ 
$\,\hskip1.4pt\omega$ on $\,M$, defined, at each point, only up to a sign, 
and having rank $\,2\,$ at every point. The image of\/ $\,\hskip1.4pt\omega\,$ 
is the Ol\-szak distribution\/ $\,\mathcal{D}$.
\end{lemma}
\begin{proof}This is obvious from  
\cite[Lemma~2.1(iii)]{derdzinski-roter:dr07c} and 
\cite[Lemma~17.1(iii)]{derdzinski-roter:dr07a} along with our characterization 
of $\,\mathcal{D}\hskip-.3pt_x$ in terms of images of nonzero $\,2$-forms.
\end{proof}
If $\,(M,g)\,$ has parallel Weyl tensor, $\,d=2$, and $\,\dim M=4$, then 
$\,\mathcal{D}$, the image of $\,\hskip1.4pt\omega\,$ (see 
Lemma~\ref{derdzinski-roter:detwo}), is a null two-di\-men\-sion\-al 
distribution, and so, by \cite[Lemma~37.8]{derdzinski-roter:dv}, at every 
$\,x\in M\,$ there exists a unique orientation of $\,T_xM\,$ for which 
$\,\omega_x$ is self-du\-al. Consequently, $\,M$ must be orientable, and a 
suitably chosen orientation makes the Weyl tensor self-du\-al.

Dunajski and West \cite{derdzinski-roter:dw} constructed various examples of 
self-du\-al metrics of the neutral signature 
$\,-\hskip.3pt-\hskip.3pt+\hskip2.4pt+\,$ in dimension four. Some of their 
metrics have parallel Weyl tensor; for instance, this is the case if $\,Q\,$ 
in \cite[formula (6.26)]{derdzinski-roter:dw} is quadratic in $\,X$.

\section{The local structure\hskip1pt: case 
$\,\,d\hskip1.2pt=1$}\label{derdzinski-roter:lsdo}

Let the data $\,I\hskip-.7pt,f,n,V\hskip-.7pt,\langle\,,\rangle,A\,$ consist of
\begin{enumerate}
  \def\theenumi{{\rm\roman{enumi}}}
\item[(i)] an open interval $\,I\subset\mathbf{R}\hskip.7pt$, a 
$\,C^\infty$ function $\,f:I\to\mathbf{R}\hskip.7pt$, and an integer $\,n\ge4$,
\item[(ii)] a real vector space $\,V\hskip.7pt$ of dimension $\,n-2\,$ 
with a pseu\-\hbox{do\hskip.7pt-}Eu\-clid\-e\-an inner product 
$\,\langle\,,\rangle$,
\item[(iii)] a nonzero traceless linear operator $\,A:V\to V\hskip-.7pt$, 
self-ad\-joint relative to $\,\langle\,,\rangle$.
\end{enumerate}
As in \cite[Theorem~3]{derdzinski-roter:r}, we use 
$\,I\hskip-.7pt,f,n,V\hskip-.7pt,\langle\,,\rangle,A\,$ to construct the 
pseu\-\hbox{do\hskip.7pt-}Riem\-ann\-i\-an manifold
\begin{equation}\label{derdzinski-roter:rcr}
(I\times\mathbf{R}\times V,\,\,\kappa\,dt^2\hskip.7pt
+\,dt\,ds\,+\,\delta)
\end{equation}
of dimension $\,n$. Here products of differentials stand for symmetric 
products, $\,t,s\,$ are the Cartesian coordinates on the 
$\,I\times\mathbf{R}\,$ factor, $\,\delta\,$ is the pull\-back to 
$\,I\times\mathbf{R}\hskip.7pt\times V\hskip.7pt$ of the flat 
pseu\-\hbox{do\hskip.7pt-}Riem\-ann\-i\-an metric on $\,V$ corresponding to 
the inner product $\,\langle\,,\rangle$, and 
$\,\kappa:I\times\mathbf{R}\hskip.7pt\times V\hskip-.7pt\to\mathbf{R}\,$ is 
the function with 
$\,\kappa(t,s,v)=f(t)\hskip.4pt\langle v,v\rangle+\langle Av,v\rangle$. 
Denoting by $\,d\,$ the dimension of the Ol\-szak distribution, we can 
characterize the manifolds (\ref{derdzinski-roter:rcr}) as follows. A proof 
can be found in \cite[Theorem~4.1]{derdzinski-roter:dr07c}:
\begin{theorem}\label{derdzinski-roter:cldeo}For any 
$\,I\hskip-.7pt,f,n,V\hskip-.7pt,\langle\,,\rangle,A\,$ as in\/ {\rm(i)} -- 
{\rm(iii)}, the pseu\-\hbox{do\hskip.7pt-}Riem\-ann\-i\-an manifold\/ 
{\rm(\ref{derdzinski-roter:rcr})} has parallel Weyl tensor and\/ $\,d=1$. In 
particular, {\rm(\ref{derdzinski-roter:rcr})} is never con\-for\-mal\-ly flat.

Conversely, in any pseu\-\hbox{do\hskip.7pt-}Riem\-ann\-i\-an manifold with 
parallel Weyl tensor and\/ $d=1$, every point has a neighborhood isometric to 
an open subset of a manifold\/ {\rm(\ref{derdzinski-roter:rcr})} constructed 
as above from some such data 
$\,I\hskip-.7pt,f,n,V\hskip-.7pt,\langle\,,\rangle,A$.

The manifold\/ {\rm(\ref{derdzinski-roter:rcr})} is locally symmetric if and 
only if\/ $\,f\,$ is constant.
\end{theorem}
One calls a pseu\-\hbox{do\hskip.7pt-}Riem\-ann\-i\-an manifold {\it 
Ric\-ci-re\-cur\-rent\/} if, for every tangent vector field $\,v$, the Ric\-ci 
tensor $\,\rho\,$ and the covariant derivative $\,\nabla_{\!v}\rho\,$ are 
linearly dependent at each point. 

Every ECS metric with $\,\,d=1\hskip1pt\,$ is Ric\-ci-re\-cur\-rent. In fact, 
according to \cite[Lemma~2.2(b)]{derdzinski-roter:dr07c}, $\,\mathcal{D}$ 
contains, at each point, the image of the Ric\-ci tensor $\,\rho$, and hence 
also the image of $\,\nabla_{\!v}\rho\,$ for any vector field $\,v\,$ (as 
$\,\mathcal{D}\,$ is parallel). Thus, $\,\rho\,$ and $\,\nabla_{\!v}\rho$, 
being symmetric, must be linearly dependent at every point if $\,d=1$. (Note 
that, in view of (\ref{derdzinski-roter:dot}) and the inclusion just 
mentioned, $\,\hskip.7pt\text{\rm rank}\hskip1.2pt\,\rho\le2\,$ at each point 
of any ECS manifold.)

Theorem~\ref{derdzinski-roter:cldeo} was inspired by the second author's 
result \cite[Theorem~3]{derdzinski-roter:r}, which is a 
gen\-er\-al-po\-si\-tion version of Theorem~\ref{derdzinski-roter:cldeo}: 
rather using the condition $\,d=1$, it assumes that the metric is 
Ric\-ci-re\-cur\-rent, $\,\rho\otimes\hskip-1pt\nabla\hskip-1pt\rho\ne0\,$ at 
every point, and $\,f\hskip2ptd\hskip-.8ptf\hskip-.6pt/\hskip-.2ptdt\ne0\,$ 
everywhere in $\,I\hskip-.7pt$.

Among ECS manifolds with $\,d=2$, in any dimension $\,n\ge4$, there are both 
Ric\-ci-re\-cur\-rent and non-Ric\-ci-re\-cur\-rent ones. (See 
\cite[Section 24]{derdzinski-roter:dr07a}; the condition 
$\,\hskip.7pt\text{\rm rank}\hskip2.7ptW\hskip-2.7pt=\hskip-1.2pt1\,$ used 
there is, by \cite[Lemma~17.1(iii)]{derdzinski-roter:dr07a}, equivalent to 
the relation $\,W\hskip-1.7pt=\pm\hskip2.5pt\omega\otimes\omega\,$ 
in Lemma~\ref{derdzinski-roter:detwo}, and hence to $\,d=2$.)

The local co\-ho\-mo\-ge\-ne\-i\-ty of any ECS metric is at most equal to 
the dimension $\,d\,$ of the Ol\-szak distribution. If $\,d=1\,$ (or, 
$\,d=2$), this follows from Theorem~\ref{derdzinski-roter:cldeo}, cf.\ 
\cite[Lemma~2.2]{derdzinski-roter:dr07b} (or, from 
Theorem~\ref{derdzinski-roter:cldet} below, cf.\ 
\cite[Re\-mark~22.1]{derdzinski-roter:dr07a} and 
\cite[the com\-ment after (f) in Section 10]{derdzinski-roter:dr07d}).

\section{The local structure\hskip1pt: case 
$\,\,d\hskip1.2pt=2$}\label{derdzinski-roter:lsdt}
Let $\,(Q,\hskip.7pt\text{\rm D}\hskip.4pt,\zeta,n,\varepsilon,V\hskip-.7pt,
\langle\,,\rangle)\,$ be a septuple formed by
\begin{enumerate}
  \def\theenumi{{\rm\alph{enumi}}}
\item[(a)] a surface $\,Q\,$ with a projectively flat tor\-sion\-free 
connection $\,\hskip.7pt{\rm D}\hskip.4pt$,
\item[(b)] a $\,\hskip.7pt\text{\rm D}\hskip.7pt$-par\-al\-lel area form 
$\,\zeta\,$ on $\,Q$, an integer $\,n\ge4$, and a sign factor 
$\,\varepsilon=\pm1$,
\item[(c)] a real vector space $\,V\hskip.7pt$ of dimension $\,n-4\,$ 
with a pseu\-\hbox{do\hskip.7pt-}Eu\-clid\-e\-an inner product 
$\,\langle\,,\rangle$.
\end{enumerate}
Also, let a twice-con\-tra\-var\-i\-ant symmetric tensor $\,\phi\,$ on $\,Q\,$ 
satisfy the differential equation
\begin{equation}\label{derdzinski-roter:ddt}
\text{\rm div}{}^{\hskip.7pt\text{\rm D}}(\text{\rm div}{}^{\hskip.7pt
\text{\rm D}}\phi)\,+\,(\hskip.4pt\rho^{\hskip.7pt\text{\rm D}}\hskip-1pt,\phi\hskip.4pt)\,\,=\,\,\varepsilon\hskip.7pt,
\end{equation}
where $\,\rho^{\hskip.7pt\text{\rm D}}$ is the Ric\-ci tensor of 
$\,\hskip.7pt\text{\rm D}\hskip.7pt$ (in coordinates: 
$\,\phi^{\hskip1.2ptjk}{}_{,\hskip.7pt jk}+\phi^{\hskip1.2ptjk}\hskip-1ptR_{jk}
\hskip.4pt=\hskip1.2pt\varepsilon$).

We define $\,\tau\,$ to be the twice-co\-var\-i\-ant symmetric tensor field 
on $\,Q\,$ corresponding to $\,\phi\,$ under the isomorphism 
$\,T\hskip-.3ptQ\to T\hskip.5pt^*\hskip-.8ptQ\,$ provided by $\,\zeta$. In 
coordinates, $\,\tau_{jk}
=\zeta_{jl}\hskip.4pt\zeta_{km}\phi^{\hskip.7pt lm}\hskip-1pt$.

Next, let $\,h^{\text{\rm D}}$ be the {\it 
Patter\-son\hskip.7pt-\hskip-.7pt Walk\-er Riemann extension metric\/} on 
$\,T\hskip.5pt^*\hskip-.8ptQ$, obtained \cite{derdzinski-roter:pw} by 
requiring that all vertical and all 
$\,\hskip.7pt\text{\rm D}\hskip.7pt$-hor\-i\-zon\-tal vectors be 
$\,h^{\text{\rm D}}\hskip-1pt$-null, while 
$\,h_x^{\text{\rm D}}(\xi,w)=\xi\hskip.7pt(d\pi_xw)\,$ for every 
$\,x\in T\hskip.5pt^*\hskip-.8ptQ$, every vector 
$\,w\in T_xT\hskip.5pt^*\hskip-.8ptQ$, and every vertical vector 
$\,\xi\in\mathrm{Ker}\hskip2.7ptd\pi_x=T_{\hskip-1pt\pi(x)}^*Q$, where 
$\,\pi:T\hskip.5pt^*\hskip-.8ptQ\to Q\,$ is the bundle projection.

Finally, we denote by $\,\delta\,$ the constant 
pseu\-\hbox{do\hskip.7pt-}Riem\-ann\-i\-an metric on $\,V$ corresponding to 
the inner product $\,\langle\,,\rangle$, and let $\,\theta\,$ stand for the 
function $\,V\hskip-.7pt\to\mathbf{R}\,$ with $\,\theta(v)=\langle v,v\rangle$.

Our septuple $\,(Q,\hskip.7pt\text{\rm D}\hskip.4pt,\zeta,n,\varepsilon,
V\hskip-.7pt,\langle\,,\rangle)\,$ now gives rise to the 
pseu\-\hbox{do\hskip.7pt-}Riem\-ann\-i\-an manifold
\begin{equation}\label{derdzinski-roter:hgt}
(T\hskip.5pt^*\hskip-.8ptQ\,\times\,V,\,\,h^{\text{\rm D}}\hskip-1.9pt
-2\tau+\hskip.7pt\delta-\theta\rho^{\hskip.7pt\text{\rm D}})\,,
\end{equation}
of dimension $\,n$, with the metric 
$\,h^{\text{\rm D}}\hskip-2.5pt-2\tau+\delta
-\theta\rho^{\hskip.7pt\text{\rm D}}\hskip-1pt$, where the function 
$\,\theta\,$ and covariant tensor fields 
$\,h^{\text{\rm D}}\hskip-1pt,\tau,\rho^{\hskip.7pt\text{\rm D}}\hskip-1pt,
\delta\,$ on $\,T\hskip.5pt^*\hskip-.8ptQ,\hskip1ptQ\,$ or $\,V\hskip-.7pt$ 
are identified with their pull\-backs to 
$\,T\hskip.5pt^*\hskip-.8ptQ\,\times\,V\hskip-1pt$.

With $\,d\,$ again denoting the dimension of the Ol\-szak distribution, we 
have the following result.
\begin{theorem}\label{derdzinski-roter:cldet}The 
pseu\-\hbox{do\hskip.7pt-}Riem\-ann\-i\-an manifold\/ 
\hskip.7pt{\rm(\ref{derdzinski-roter:hgt})} obtained as above from any given 
septuple $\,(Q,\hskip.7pt\text{\rm D}\hskip.4pt,\zeta$, 
$\,n,\varepsilon,V\hskip-.7pt,\langle\,,\rangle)\,$ with\/ {\rm(a)} -- 
{\rm(c)} has parallel Weyl tensor and\/ $\,d=2$.

Conversely, in any pseu\-\hbox{do\hskip.7pt-}Riem\-ann\-i\-an manifold with 
parallel Weyl tensor and\/ $d=2$, every point has a neighborhood isometric to 
an open subset of a manifold\/ \hskip.7pt{\rm(\ref{derdzinski-roter:hgt})} 
constructed from some septuple\/ 
$\,(Q,\hskip.7pt\text{\rm D}\hskip.4pt,\zeta,n,\varepsilon,V\hskip-.7pt
,\langle\,,\rangle)\,$ as in\/ {\rm(a)} -- {\rm(c)}.

The manifold\/ \hskip.7pt{\rm(\ref{derdzinski-roter:hgt})} is never con\-for\-mal\-ly flat, and 
it is locally symmetric if and only if the Ric\-ci tensor\/ 
$\,\rho^{\hskip.7pt\text{\rm D}}$ is 
$\,\hskip.7pt\text{\rm D}\hskip.7pt$-par\-al\-lel.
\end{theorem}
Our septuple $\,(Q,\hskip.7pt\text{\rm D}\hskip.4pt,\zeta,n,\varepsilon,
V\hskip-.7pt,\langle\,,\rangle)\,$ of parameters does not include $\,\phi$, 
even though the metric in (\ref{derdzinski-roter:hgt}) evidently depends on 
$\,\hskip1pt\tau\,$ (and hence on $\,\hskip.7pt\phi$). This is justified by 
the fact that, locally, $\,\phi$ with (\ref{derdzinski-roter:ddt}) always 
exists, and, for any fixed 
$\,(Q,\hskip.7pt\text{\rm D}\hskip.4pt,\zeta,n,\varepsilon,V\hskip-.7pt,
\langle\,,\rangle)$, the manifolds (\ref{derdzinski-roter:hgt}) corresponding 
to two choices of $\,\phi\,$ are, locally, isometric to each other 
\cite[Re\-mark~22.1]{derdzinski-roter:dr07a}.

In dimension four, $\,V\hskip-.7pt=\{0\}\,$ and $\,\langle\,,\rangle=0$, so 
that the septuple $\,(Q,\hskip.7pt\text{\rm D}\hskip.4pt,\zeta,n,\varepsilon,
V\hskip-.7pt,\langle\,,\rangle)\,$ with (a) -- (c) may be replaced by a 
quadruple $\,(Q,\hskip.7pt\text{\rm D}\hskip.4pt,\zeta,\varepsilon)\,$ 
consisting of a surface $\,Q$, a projectively flat tor\-sion\-free connection 
$\,\hskip.7pt{\rm D}\hskip.7pt\,$ on $\,Q$, a 
$\,\hskip.7pt\text{\rm D}\hskip.7pt$-par\-al\-lel area form $\,\zeta\,$ on 
$\,Q$, and a sign factor $\,\varepsilon=\pm1$. The pair 
(\ref{derdzinski-roter:hgt}) then becomes the 
pseu\-\hbox{do\hskip.7pt-}Riem\-ann\-i\-an four-man\-i\-fold 
$\,(T\hskip.5pt^*\hskip-.8ptQ\hskip1pt,\,\,
h^{\text{\rm D}}\hskip-1.9pt-2\tau)$.

\section{Compact ECS manifolds\hskip1pt: 
existence}\label{derdzinski-roter:cxis}
\begin{theorem}\label{derdzinski-roter:exist}In every dimension $\,n\ge5\,$ 
with 
\hbox{$n$\hskip3pt$\equiv$\hskip4pt$5$\hskip4.4pt$(${\rm 
mod}\hskip2.7pt$3\hskip.5pt)$}, there exists a compact ECS manifold of any 
prescribed indefinite metric signature, dif\-feo\-mor\-phic to a nontrivial 
torus bundle over the circle.
\end{theorem}
A proof of Theorem~\ref{derdzinski-roter:exist} is outlined in 
\S\ref{derdzinski-roter:pfth}. For a more detailed argument, see 
\cite{derdzinski-roter:dr07b}.

The compact ECS manifolds that are shown to exist in 
\S\ref{derdzinski-roter:pfth} have further properties, not included in the 
statement of Theorem~\ref{derdzinski-roter:exist}. Some of these properties 
are mentioned in \S\ref{derdzinski-roter:comm} (items II, IV and V). The first 
four steps of the proof in \S\ref{derdzinski-roter:pfth} might {\it a priori} 
produce ECS metrics on bundles over the circle, the fibre of which is either a 
torus, or a $\,2$-step nil\-man\-i\-fold admitting a complete flat 
tor\-sion\-free connection with a nonzero parallel vector field. However, our 
particular existence argument realizes only the torus as the fibre. See 
\cite[Re\-marks~4.1(iv) and~6.2]{derdzinski-roter:dr07b}.

\section{Compact ECS manifolds\hskip1pt: 
properties}\label{derdzinski-roter:cprp}
First, the existence of an ECS metric on a given compact manifold imposes 
specific restrictions on the fundamental group, Euler characteristic, and real 
Pontryagin classes:
\begin{theorem}\label{derdzinski-roter:pichp}If a compact manifold\/ $\,M\,$ 
of dimension $\,n\ge4\,$ admits an ECS metric, then\/ $\,\pi_1M\,$ is 
infinite, $\,\chi(M)=0$, and\/ $\,p_i(M)=0\,$ in $\,H^{4i}(M,\mathbf{R})\,$ 
for all\/ $\,i\ge1$.
\end{theorem}
Here $\,\chi(M)\,$ and $\,p_i(M)\,$ vanish since so do the 
Gauss-Bon\-net-Chern integrand and the Pont\-rya\-gin forms 
\cite{derdzinski-roter:dr07d}. For proofs of infiniteness of 
$\pi_1M\hskip-1pt$ and the next two theorems, see 
\cite{derdzinski-roter:dr07d} as well.
\begin{theorem}\label{derdzinski-roter:lornc}Every four-di\-men\-sion\-al 
Lo\-rentz\-i\-an ECS manifold is noncompact.
\end{theorem}
\begin{theorem}\label{derdzinski-roter:lorbc}Let\/ $\,(M,g)\,$ be a compact 
Lo\-rentz\-i\-an ECS manifold. Then some two-fold covering manifold of\/ 
$\,M\hskip.7pt$ is the total space of a $\,C^\infty\hskip-1pt$ bundle over the 
circle, the fibre of which admits a flat tor\-sion\-free connection with a 
nonzero parallel vector field.
\end{theorem}
If $\,(M,g)\,$ is a compact ECS manifold and $\,d\in\{1,2\}\,$ is the 
dimension of its Ol\-szak distribution (\S\ref{derdzinski-roter:odis}), then 
$\,T\hskip-.3ptM=H\oplus H^+\hskip-1.6pt\oplus H^-$ for some vector 
sub\-bun\-dles $\,H,H^+$ and $\,H^-$ of $\,T\hskip-.3ptM\,$ such that $\,H^+$ 
is space\-like, $\,H^-$ is time\-like, and  both $\,H^\pm$ have the fibre 
dimension $\,d$.

In fact, we obtain $\,H^+$ and $\,H^-$ by decomposing $\,T\hskip-.3ptM\,$ into 
a space\-like and a time\-like sub\-bundle, and then projecting the null 
distribution $\,\mathcal{D}\,$ onto the summands.

\section{Some open questions}\label{derdzinski-roter:open}
\begin{enumerate}
  \def\theenumi{{\rm\Roman{enumi}}}
\item Do compact ECS manifolds exist in dimension four? 
\item Does any torus admit an ECS metric? More generally, does 
there exist a compact ECS manifold with an Abelian fundamental group?
\item Are there compact ECS manifolds of dimensions $\,n\ge5\,$ other than 
$\,n=3j+2,\,j\in\mathbf{Z}\hskip1pt$?
\item Can a compact ECS manifold be locally homogeneous?
\item Must the Ol\-szak distribution of a compact ECS manifold be 
one-di\-men\-sion\-al? More generally, are all compact ECS manifolds 
Ric\-ci-re\-cur\-rent?
\end{enumerate}

\section{Comments on Questions I -- V in 
\S\ref{derdzinski-roter:open}}\label{derdzinski-roter:comm}
\begin{enumerate}
  \def\theenumi{{\rm\Roman{enumi}}}
\item If the answer to Question I is `yes' and compact four-di\-men\-sion\-al 
ECS manifolds do exist, they all must have the neutral metric signature 
$\,-\hskip.3pt-\hskip.3pt+\hskip2.4pt+\hskip.7pt$. In fact, they can be 
neither Riemannian \cite[Theorem~2]{derdzinski-roter:dr77}, nor 
Lo\-rentz\-i\-an (Theorem~\ref{derdzinski-roter:lornc}). See also item V below.
\item None of the known compact ECS manifolds admits a finite covering by a 
manifold with an Abelian fundamental group.
\item The matrices (\ref{derdzinski-roter:mtr}) acting in 
$\,\mathbf{R}\hskip-.7pt^3$ lead to a translation operator that has the 
property required in STEP~\ref{derdzinski-roter:lttce} of 
\S\ref{derdzinski-roter:pfth}. The $\,j$th Car\-te\-sian-pow\-er extension of 
(\ref{derdzinski-roter:mtr}) does the same in 
$\,\mathbf{R}\hskip-.7pt^{3j}\hskip-1pt$, which is why our argument yields 
compact ECS manifolds of dimensions $\,n=3j+2$. A ``building block'' 
$\,\mathbf{R}\hskip-.7pt^m\hskip-1pt$, $\,m\ge4$, instead of 
$\,\mathbf{R}\hskip-.7pt^3\hskip-1pt$, with operators analogous to 
(\ref{derdzinski-roter:mtr}), would lead to examples in other dimensions, and 
it seems reasonable to expect such operators to exist, although they might be 
harder to find than those appearing in (\ref{derdzinski-roter:mtr}).
\item Locally homogeneous ECS manifolds exist \cite{derdzinski-roter:d} in 
every dimension $\,n\ge4$. However, none of the compact ECS manifolds arising 
from the argument in \S\ref{derdzinski-roter:pfth} is locally homogeneous.
\item Every known compact ECS manifold has a one-di\-men\-sion\-al Ol\-szak 
distribution, and is thereforee Ric\-ci-re\-cur\-rent 
(\S\ref{derdzinski-roter:lsdo}). One might try to answer Question I 
in the affirmative and, simultaneously, the first part of Question V in the 
negative, by proceeding as follows. We begin by choosing a quadruple 
$\,(Q,\hskip.7pt\text{\rm D}\hskip.4pt,\zeta,\varepsilon)\,$ as at the end of 
\S\ref{derdzinski-roter:lsdt}, with a {\it closed\/} 
surface $\,Q$. (According to \cite[Section 23]{derdzinski-roter:dr07a}, such 
quadruples exist, and realize all dif\-feo\-mor\-phic types of closed surfaces 
$\,Q$.) As a next step, we need to find a discrete group $\,\Gamma\hskip.4pt$ 
of isometries of the four-di\-men\-sion\-al ECS manifold 
$\,(T\hskip.5pt^*\hskip-.8ptQ\hskip1pt,\,\,
h^{\text{\rm D}}\hskip-1.9pt-2\tau)$, acting on 
$\,T\hskip.5pt^*\hskip-.8ptQ\,$ properly dis\-con\-tin\-u\-ous\-ly with a 
compact quotient. There is a topological restriction: $\,\Gamma\hskip.4pt$ as 
above cannot exist unless $\,Q\,$ is dif\-feo\-mor\-phic to the torus $\,T^2$ 
or the Klein bottle $\,K^2$ (see below). We do not know if such 
$\,\Gamma\hskip.4pt$ exists either for $\,T^2$ or for $\,K^2\hskip-1pt$, even 
though, on $\,T^2\hskip-1pt$, the connections in question are relatively well 
understood, cf.\ 
\cite[the dis\-cus\-sion fol\-low\-ing Re\-mark~7.2]{derdzinski-roter:dr07a}.
\end{enumerate}
The ``topological restriction'' mentioned under V can be phrased as 
follows.

If $\,Q\,$ is a compact manifold and some group $\,\Gamma\hskip.4pt$ of 
dif\-feo\-mor\-phisms of $\,T\hskip.5pt^*\hskip-.8ptQ\,$ acts on 
$\,T\hskip.5pt^*\hskip-.8ptQ\,$ properly dis\-con\-tin\-u\-ous\-ly with a 
compact quotient, then $\,\chi(Q)=0$. In fact, we may assume that $\,Q\,$ is 
orientable (by passing, if necessary, to a two-fold orientable covering and 
replacing $\,\Gamma\hskip.4pt$ with a $\,\mathbf{Z}_2$ extension). Since the inclusion 
of the zero section $\,Q\,$ in $\,T\hskip.5pt^*\hskip-.8ptQ\,$ is a homotopy 
equivalence, a generator $\,[Q]\,$ of 
$\,H_m(Q,\mathbf{Z})\,\approx\,\mathbf{Z}$, for $\,m=\dim Q$, is also a 
generator of $\,H_m(T\hskip.5pt^*\hskip-.8ptQ,\mathbf{Z})$, and so its image 
$\,F_*[Q]\,$ under any dif\-feo\-mor\-phism 
$\,F:T\hskip.5pt^*\hskip-.8ptQ\to T\hskip.5pt^*\hskip-.8ptQ\hskip1pt\,$ equals 
$\,\hskip.7pt\pm\hskip.7pt[Q]$. We now choose $\,F\in\Gamma$ such that $\,Q\,$ 
and $\,F(Q)\,$ are disjoint: due to compactness of $\,Q\,$ and proper 
discontinuity of the action, this is the case for all but finitely many 
$\,F\in\Gamma$. Thus, $\,\chi(Q)=[Q]\cdot[Q]=\pm\hskip.7pt[Q]\cdot F_*[Q]=0$, 
where $\,\cdot\,$ is the intersection form in 
$\,H_m(T\hskip.5pt^*\hskip-.8ptQ,\mathbf{Z})$. (The same argument remains 
valid if $\,T\hskip.5pt^*\hskip-.8ptQ$ is replaced by the total space of an 
orientable real vector bundle of fibre dimension $\,n\,$ over a 
compact orientable $\,n$-di\-men\-sion\-al manifold $\,Q$, the conclusion 
being now that the Euler number of the vector bundle must vanish.)

\section{Proof of Theorem~\ref{derdzinski-roter:exist} (an outline)}\label{derdzinski-roter:pfth}
\begin{step}\label{derdzinski-roter:group}Consider an $\,n$-di\-men\-sion\-al 
pseu\-\hbox{do\hskip.7pt-}Riem\-ann\-i\-an ECS manifold 
(\ref{derdzinski-roter:rcr}) constructed from some 
$\,I\hskip-.7pt,f,n,V\hskip-.7pt,\langle\,,\rangle,A\,$ with (i) -- (iii) in 
\S\ref{derdzinski-roter:lsdo} such that 
$\,I\hskip-.7pt=\mathbf{R}\hskip.7pt$, while $\,f\,$ is nonconstant (cf.\ 
Theorem~\ref{derdzinski-roter:cldeo}) and periodic, for some period $\,p>0$. 
Let $\,\mathcal{E}$ be the vector space of all $\,C^\infty$ solutions 
$\,u:\mathbf{R}\to V\hskip.4pt$ to the differential equation 
$\,\ddot u(t)=f(t)\hskip.4pt u(t)+Au(t)$. The set 
$\,\mathrm{G}=\mathbf{Z}\times\mathbf{R}\times\mathcal{E}\,$ has a group 
structure such that, setting $\,(k,q,u)\cdot(t,s,v)
=(t+kp,\hskip.7pt s+q-\langle\dot u(t),2v+u(t)\rangle,\hskip.7pt v+u(t))$, for 
$\,(k,q,u)\in\hskip.4pt\mathrm{G}\hskip.7pt$ and 
$\,(t,s,v)\in\hskip.4pt\mathbf{R}\hskip-.5pt^2\hskip.7pt\times\,V\hskip-.7pt$, 
we define an action of $\,\mathrm{G}\hskip.7pt$ on 
$\,\mathbf{R}\hskip-.5pt^2\hskip.7pt\times\,V$ from the left. The 
action of $\,\mathrm{G}\hskip.7pt$ consists of isometries of our manifold 
(\ref{derdzinski-roter:rcr}). See \cite[Lemma~2.2]{derdzinski-roter:dr07b}.
\end{step}
\begin{step}\label{derdzinski-roter:prdis}Question: when does a manifold 
(\ref{derdzinski-roter:rcr}) selected as in STEP~\ref{derdzinski-roter:group}, 
with the corresponding 
$\,f,n,V\hskip-.7pt,\langle\,,\rangle,A,p,\mathcal{E},\hskip.7pt\mathrm{G}\,$ 
and $\,I\hskip-.7pt=\mathrm{R}\hskip.7pt$, lead to a compact ECS manifold that 
arises as the quotient of (\ref{derdzinski-roter:rcr}) under a discrete 
subgroup $\,\Gamma\hskip.7pt$ of $\,\mathrm{G}\hskip.7pt$ acting on 
$\,\mathbf{R}\hskip-.5pt^2\hskip.7pt\times\,V$ properly 
dis\-con\-tin\-u\-ous\-ly? Answer: such $\,\Gamma\hskip.7pt$ exists if and 
only if some $\,C^\infty$ curve 
$\,\mathbf{R}\ni t\mapsto B(t)\in\hskip.7pt\text{\rm End}\hskip.7pt(V)$, 
periodic of period $\,p$, satisfies the differential equation 
$\,\dot B(t)\hskip.7pt+[B(t)]^2\hskip-.7pt=f(t)\hskip.7pt+A\,$ (briefly, 
$\,\dot B+B^2\hskip-.7pt=f\hskip.7pt+A$) and, at the same time, certain 
additional conditions, listed in 
\cite[Theorem~6.1(ii)]{derdzinski-roter:dr07b}, hold for some lattice 
$\,\Sigma\,$ in the vector space $\,\mathcal{W}=\mathbf{R}\times\mathcal{L}$, 
where $\,\mathcal{L}\,$ is the space of all $\,C^\infty$ functions 
$\,u:\mathbf{R}\to V$ with $\,\dot u(t)=B(t)\hskip.4ptu(t)$. (Thus, 
$\,\mathcal{L}\subset\mathcal{E}$.) The most important of these additional 
conditions is an arithmetic property of the {\it translation operator\/} 
$\,T:\mathcal{L}\to\mathcal{L}\,$ acting by $\,(Tu)(t)=u(t-p)$. Namely, we 
require the existence of a linear functional $\,\varphi\in\mathcal{L}^*$ such 
that $\,\varPsi(\Sigma)=\Sigma\,$ for the operator 
$\,\varPsi:\mathcal{W}\to\mathcal{W}\,$ given by 
$\,\varPsi(r,u)=(r+\varphi(u),Tu)$.
\end{step}
\begin{step}\label{derdzinski-roter:lttce}Observe that, in 
STEP~\ref{derdzinski-roter:prdis}, a lattice $\,\Sigma\,$ with 
$\,\varPsi(\Sigma)=\Sigma\,$ for some functional $\,\varphi\,$ exists if and 
only if $\,\text{\rm det}\,T=\pm1\,$ and the matrix of $\,T\,$ in some basis 
of $\,\mathcal{L}\,$ consists of integers, cf.\ 
\cite[the end of Section~1]{derdzinski-roter:dr07b}. Thus, the main part of 
our task is to find $\,B\,$ with $\,\dot B+B^2\hskip-.7pt=f\hskip.7pt+A\,$ 
such that the corresponding translation operator $\,T\,$ has the property just 
stated.
\end{step}
\begin{step}\label{derdzinski-roter:fixkl}Fix integers $\,k\,$ and $\,l\,$ 
with $\,4\,<\,k\,<\,l\,\le\,k^2\hskip-1pt/4$. (Such $\,l\,$ exists for any 
given $\,k>4$.) It is an easy exercise 
\cite[Lemma~1.3]{derdzinski-roter:dr07b} to verify that the polynomial 
$\,\hskip1.5ptP(\lambda)\,=\,-\lambda^3\hskip.7pt
+\,k\hskip.7pt\lambda^2\hskip.7pt-\,l\hskip.7pt\lambda\,+\,1$ then has real 
roots $\,\lambda,\mu,\nu\,$ with $\,0<\lambda<\mu<\nu$,\hskip6pt$\lambda<1<\nu$,\hskip6pt$\lambda\mu<1<\mu\nu\,$ 
and $\,\lambda\nu\ne1$.
\end{step}
\begin{step}\label{derdzinski-roter:fixpp}Fix $\,p\in(0,\infty)\,$ and denote 
by $\,\mathcal{F}_{\hskip-.7ptp}$ the set of all septuples 
$\,(\alpha,\beta,\gamma,f,a,b,c)\,$ consisting of $\,C^\infty$ functions 
$\,\alpha,\beta,\gamma,f:\mathbf{R}\to\mathbf{R}\,$ of the variable $\,t$, 
periodic of period $\,p$, and three distinct real constants $\,a,b,c\,$ with 
$\,a+b+c=0$, of which $\,b\,$ is the smallest, such that, with 
$\,(\,\,)\dot{\,}=\,d/dt$,
\begin{equation}\label{derdzinski-roter:thr}
\dot\alpha+\alpha^2=f+\hskip.7pta\hskip.7pt,\hskip16pt\dot\beta+\beta^2
=f+\hskip.7ptb\hskip.7pt,\hskip16pt
\dot\gamma+\gamma^2=f+\hskip.7ptc\hskip.7pt,
\end{equation}
and $\,\alpha>\beta>\gamma$. Let $\,\mathcal{C}\,$ be the subset of 
$\,\mathcal{F}_{\hskip-.7ptp}$ formed by all 
$\,(\alpha,\beta,\gamma,f,a,b,c)\,$ with {\it constant\/} 
$\,\alpha,\beta,\gamma\,$ and $\,f$. Define a mapping 
$\,\hskip.7pt\text{\rm spec}:\mathcal{F}_{\hskip-.7ptp}
\to\mathbf{R}\hskip-.7pt^3$ by 
$\,\hskip.7pt\text{\rm spec}\hskip.7pt(\alpha,\beta,\gamma,f,a,b,c)\,
=\,(\lambda,\mu,\nu)$, where
\begin{equation}\label{derdzinski-roter:exp}
(\lambda,\mu,\nu)\,\,
=\,\,(\exp\,[-\textstyle{\int_0^p\alpha(t)\,dt}\hskip.7pt]\hskip.7pt,
\hskip.7pt\,\,\exp\,[-\textstyle{\int_0^p\beta(t)\,dt}\hskip.7pt]\hskip.7pt,
\hskip.7pt\,\,\exp\,[-\textstyle{\int_0^p\gamma(t)\,dt}\hskip.7pt]
\hskip.7pt)\hskip.7pt.
\end{equation}
\end{step}
\begin{step}\label{derdzinski-roter:showt}Show that 
$\,\hskip.7pt\text{\rm spec}\hskip.7pt(\mathcal{F}_{\hskip-.7ptp}\hskip-1pt
\smallsetminus\mathcal{C})\,=\,\mathcal{U}\hskip-.4pt$, where 
$\,\mathcal{U}\subset\mathbf{R}\hskip-.7pt^3$ is the open set of all 
$\,(\lambda,\mu,\nu)\,$ with $\,0<\lambda<\mu<\nu$,\hskip6pt$\lambda<1<\nu$,\hskip6pt$\lambda\mu<1<\mu\nu\,$ 
and $\,\lambda\nu\ne1$. This is done by 
reducing the number of unknown functions in a septuple with 
(\ref{derdzinski-roter:thr}). Specifically, to solve the two-e\-qua\-tion 
system $\,\dot\alpha+\alpha^2=\hskip.7ptf+\hskip.7pta$, 
$\,\dot\beta+\beta^2=\hskip.7ptf+\hskip.7ptb\,$ with $\,\alpha>\beta\,$ 
(treating our $\,a,b\,$ as fixed), one sets $\,\rho=\alpha-\beta$ and 
$\,\psi=\alpha+\beta$. As $\,\psi=(a-b-\dot\rho)/\hskip-1.2pt\rho\hskip.7pt$, 
we may reconstruct $\,\alpha\,$ and $\,\beta\,$ from $\,\rho$. Solutions 
$\,(\alpha,\beta)$ of the two-e\-qua\-tion system are thus in a bijective 
correspondence with arbitrary $\,C^\infty$ functions 
$\,\rho:\mathbf{R}\to(0,\infty)$, periodic of period $\,p$. Similarly, 
a solution to the last two equations in (\ref{derdzinski-roter:thr}) is 
represented by a single arbitrary positive function analogous to $\,\rho$, 
which we denote by $\,\sigma$. Our $\,\rho\,$ and $\,\sigma\,$ are subject to 
a single differential equation, stating that $\,\beta\,$ reconstructed from 
$\,\rho\,$ is the same as $\,\beta\,$ reconstructed from $\,\sigma$. The 
function $\,\log\hskip1.2pt
(\hskip-.4pt\sigma\hskip-.6pt/\hskip-1.1pt\rho\hskip-.2pt)\,$ may, however, be 
completely arbitrary (aside from being of class $\,C^\infty$ and periodic of 
period $\,p$); see \cite[Lemma~9.6]{derdzinski-roter:dr07b}. Expressing the 
triple (\ref{derdzinski-roter:exp}) in terms of the new unknown function 
$\,\log\hskip1.2pt(\hskip-.4pt\sigma\hskip-.6pt/\hskip-1.1pt\rho\hskip-.2pt)$, 
we now verify that 
$\,\hskip.7pt\text{\rm spec}\hskip.7pt(\mathcal{F}_{\hskip-.7ptp}\hskip-1pt
\smallsetminus\mathcal{C})\,=\,\mathcal{U}\hskip-.4pt$.
\end{step}
\begin{step}\label{derdzinski-roter:imspc}Using $\,\,(\lambda,\mu,\nu)\,\,$ 
obtained in STEP~\ref{derdzinski-roter:fixkl}, apply 
STEP~\ref{derdzinski-roter:showt} to pick 
$\,(\alpha,\beta,\gamma,f,a,b,c)\in\mathcal{F}_{\hskip-.7ptp}\hskip-1pt
\smallsetminus\mathcal{C}$ with 
$\,\hskip.7pt\text{\rm spec}\hskip.7pt(\alpha,\beta,\gamma,f,a,b,c)\,
=\,(\lambda,\mu,\nu)$, and then set
\begin{equation}\label{derdzinski-roter:mtr}
B(t)\,\,=\,\left[\begin{matrix}\alpha(t)&0&0\cr
0&\beta(t)&0\cr
0&0&\gamma(t)\end{matrix}\right],\hskip25pt
A\,\,=\,\left[\begin{matrix}a&0&0\cr
0&b&0\cr
0&0&c\end{matrix}\right].
\end{equation}
As a consequence of (\ref{derdzinski-roter:thr}), we now have 
$\,\dot B+B^2\hskip-.7pt=f\hskip.7pt+A$. Treating the matrices 
(\ref{derdzinski-roter:mtr}) as endomorphisms of the space 
$\,\mathbf{R}\hskip-.7pt^3$ endowed with the standard 
pseu\-\hbox{do\hskip.7pt-}Eu\-clid\-e\-an inner product (of any signature), 
and noting that all $\,B(t)\,$ commute with one another, we easily conclude 
that the spectrum of the corresponding translation operator 
(STEP~\ref{derdzinski-roter:prdis}) is given by (\ref{derdzinski-roter:exp}), 
and hence coincides with our fixed $\,(\lambda,\mu,\nu)$. Thus, 
$\,P(\lambda)\,$ appearing in STEP~\ref{derdzinski-roter:fixkl} is the 
characteristic polynomial of $\,T$, and, consequently, $\,T\,$ has the 
property required in STEP~\ref{derdzinski-roter:lttce} (see 
\cite[the end of Section~1]{derdzinski-roter:dr07b}). Also, $\,f\,$ is 
nonconstant since $\,(\alpha,\beta,\gamma,f,a,b,c)\notin\mathcal{C}$, cf.\ 
\cite[Re\-mark~9.1]{derdzinski-roter:dr07b}. The same conclusions hold if the 
endomorphisms (\ref{derdzinski-roter:mtr}) are replaced by their $\,j$th 
Cartesian powers acting in 
$\,V\hskip-.7pt=\mathbf{R}\hskip-.7pt^{3j}\hskip-1pt$, $\,j\ge1$. According to 
STEP~\ref{derdzinski-roter:lttce}, this completes the proof.
\end{step}

\end{document}